\documentclass[10pt]{article}
\usepackage[margin=0.9in]{geometry}
\usepackage[T1]{fontenc}
\usepackage{lmodern,amsmath,amssymb,amsthm,booktabs,array,microtype}
\usepackage[hidelinks]{hyperref}
\newtheorem{theorem}{Theorem}[section]
\newtheorem{lemma}[theorem]{Lemma}
\newtheorem{proposition}[theorem]{Proposition}
\newtheorem{corollary}[theorem]{Corollary}
\theoremstyle{definition}
\newcommand{\norm}[1]{\left\lVert#1\right\rVert_2}
\newcommand{\len}{\operatorname{len}}
\newcommand{\MST}{\operatorname{MST}}
\newcommand{\MSTthree}{\operatorname{MST}_3}
\title{Degree-three Euclidean spanning trees\\within a factor of $1.514$ of the minimum}
\author{\begin{tabular}{c}
William Han$^{1}$\quad Dan Ismailescu$^{2}$\thanks{Corresponding author: \href{mailto:dan.p.ismailescu@hofstra.edu}{dan.p.ismailescu@hofstra.edu}.}\quad Ganghun Kim$^{3}$\\[0.7em]
\small $^{1}$Horace Mann School, NY, USA\\
\small \href{mailto:william_han@horacemann.org}{\texttt{william\_han@horacemann.org}}\\[0.3em]
\small $^{2}$Mathematics Department, Hofstra University, NY, USA\\
\small \href{mailto:dan.p.ismailescu@hofstra.edu}{\texttt{dan.p.ismailescu@hofstra.edu}}\\[0.3em]
\small $^{3}$Fort Lee High School\\
\small \href{mailto:jasonghkim@gmail.com}{\texttt{jasonghkim@gmail.com}}
\end{tabular}}
\date{}
\begin{document}\maketitle
\begin{abstract}
For every finite point set $X$ in a Euclidean space of arbitrary dimension,
we prove that there is a spanning tree of maximum degree three and total
length at most $(757/500)\MST(X)=1.514\MST(X)$, where $\MST(X)$ denotes
the unrestricted minimum spanning-tree length. More strongly, for any
root and any finite list of other points, there is a Hamiltonian path
starting at the root whose length is at most $1.514$ times the sum of the
root-to-point distances. The proof follows the rooted-path and
variable-block framework used in earlier work, and applies Young's
weighted-distance inequality to explicit rational mixtures of paths.
Eight ordered radius regions cover the normalized domain. Their local
estimates require only 29 vertex comparisons, all certified by rational
arithmetic. The cover uses the three thresholds $1/4$, $3/8$, and $1/2$;
its largest block has nine vertices. We give complete proofs of the
normalization transfer, the global induction, and the terminal cases
with zero or one actual tail point. We also give a fixed-dimensional
cross-polytope lower bound, present Zbarsky's orthogonal simplex-cluster
construction in explicit coordinates, and compare their finite-dimension
bounds. The base cases complete the upper-bound argument before the
lower-bound constructions are presented.
\end{abstract}

\section{Introduction}
A minimum spanning tree connects a finite set of points with the least
possible total edge length. A degree constraint limits the number of
edges incident with a vertex and may force a longer network. For a finite
set $X$ in a Euclidean space, write $\MST(X)$ for the minimum length of
an unrestricted spanning tree and $\MSTthree(X)$ for the minimum length
among spanning trees of maximum degree at most three. Edges are straight
segments, their weights are Euclidean lengths, and no additional vertices
are introduced at crossings.

The dimension-free problem asks how large the ratio
$\MSTthree(X)/\MST(X)$ can be. A doubled-tree traversal followed by
shortcutting gives an elementary factor-two guarantee, even for a
spanning path in an arbitrary metric space. Improving this guarantee
uses Euclidean information that the triangle inequality alone does not
capture.

\subsection{History and previous work}
Papadimitriou and Vazirani~\cite{PV84} studied the computational complexity
of Euclidean degree-constrained spanning trees in 1984, including the
NP-hardness of finding a shortest degree-three tree. Their work also
motivated approximation guarantees that exploit geometry.
Khuller, Raghavachari, and Young~\cite{KRY96} obtained an
arbitrary-dimensional factor $5/3$, and the stronger factor $3/2$ in the
plane. Their results appeared in conference form in 1994 and in journal
form in 1996. An important ingredient is the replacement of child stars
in a rooted spanning tree by short paths.

Chan~\cite{Chan04} improved the arbitrary-dimensional degree-three
factor to $2\sqrt6/3\approx1.632993$, and obtained a planar factor below
$1.402$. His conference and journal versions appeared in 2003 and 2004.
The planar and dimension-free problems must be distinguished: a bound
proved using planar geometry does not automatically extend to all
Euclidean dimensions.

Zbarsky~\cite{Zbarsky14} proved the arbitrary-dimensional factor $1.559$.
His proof controls paths with a prescribed starting point, keeps two
possible ending points, and varies the size of a removed block according
to the ordered radii. He also constructed examples giving a tree-ratio
lower bound of about $1.447$. Young's weighted-distance inequality,
from a 1994 technical report~\cite{Young94}, bounds weighted sums of
pairwise distances when the point radii are prescribed. That inequality
is used in Zbarsky's proof and is also the geometric tool used here.
We prove the required upper-bound form directly.

\begin{center}\begin{tabular}{lcc}\toprule
Work&Planar degree-three factor&Dimension-free factor\\\midrule
Khuller--Raghavachari--Young~\cite{KRY96}&$3/2$&$5/3$\\
Chan~\cite{Chan04}&$<1.402$&$2\sqrt6/3$\\
Zbarsky~\cite{Zbarsky14}&Dimension-free bound applies&$1.559$\\
This paper&Dimension-free bound applies&$757/500=1.514$\\\bottomrule
\end{tabular}\end{center}
This table records the results relevant to the present argument; it is
not intended as a survey of all variants of degree-constrained network
design.

\subsection{Results}
To our knowledge, the best previously published dimension-free upper
bound was $1.559$, due to Zbarsky~\cite{Zbarsky14}. We improve this bound
to $757/500=1.514$, refining the rooted-path framework with explicit
weighted path averages and an eight-region verification.

Our exact constant throughout is
\[
C=757/500=1.514.
\]
\begin{theorem}\label{thm:tree}
For every finite set $X$ in any Euclidean space,
\[
\MSTthree(X)\le C\MST(X).
\]
\end{theorem}
For a singleton both lengths are zero. For sets with at least two points,
the theorem bounds the ratio of the two positive optimum lengths.
Equivalently, if $\tau_3^D$ is the supremum of that ratio over finite
sets in $\mathbb R^D$, then $\sup_{D\ge1}\tau_3^D\le757/500$.

The stronger path statement below implies Theorem~\ref{thm:tree}.
We state it for labeled configurations, allowing coincident locations.
This convention is useful when radii are contracted or copies of the
origin are introduced temporarily.
\begin{theorem}[Two prescribed endpoints]\label{thm:paths}
Let $a_1,\ldots,a_m$ be labeled points, $m\ge2$, and put
$d_i=\norm{a_i}$. Order the labels so that $d_1\ge\cdots\ge d_m$.
There are two paths $P_1,P_2$, each starting at the origin and visiting
all $m$ labels exactly once, ending at $a_1,a_2$, respectively, such that
\[
\frac{\len(P_1)+\len(P_2)}{2}\le C\sum_{i=1}^m d_i.
\]
\end{theorem}
A path here means a vertex order with its Euclidean edge lengths; its
straight segments may cross. The origin is a designated root label,
distinct from the listed labels even if their locations coincide.
Taking the shorter of $P_1,P_2$ proves the following consequence.
\begin{corollary}\label{cor:rooted}
For any root $v$ and any finite labeled list $a_1,\ldots,a_m$, there is
a Hamiltonian path starting at $v$ of length at most
$C\sum_i\norm{a_i-v}$.
\end{corollary}
For zero or one listed point the corollary follows directly, and for at
least two points it follows from Theorem~\ref{thm:paths} after translation.

\subsection{The new ingredients and the organization of the proof}
The rooted-path reduction follows Khuller, Raghavachari, and
Young~\cite[Section~4]{KRY96}; see also Chan~\cite[Sections~2 and~5]{Chan04}.
The two-prescribed-endpoint formulation and variable-block induction
follow Zbarsky~\cite[Section~4]{Zbarsky14}. The weighted-distance
inequality is due to Young~\cite[Lemma~1]{Young94}. The convexity
principle underlying vertex checks is standard. Our contribution is an explicit choice and
organization of the local estimates that proves the smaller constant.

First, rational mixtures of three path templates and positive rational
witness vectors retain the unequal block radii in one quadratic bound.
Second, the two incoming points are averaged together. Ordering may bound
one incoming radius by a block radius and the other by a fixed threshold;
retaining that distinction enlarges the applicable regions.
Third, a quadratic majorant handles a capped incoming-radius term in the
five-point block. After completing a square, the bound becomes the norm
of a single affine vector on the whole region. Finally, eight ordered
regions cover all normalized profiles, using only 29 vertex checks.
The required block sizes are $3,4,5,6,7,9$.

The lower-bound comparison in Section~\ref{sec:lower} explains a
different dimensional effect. The cross-polytope family gives the exact
ratio $[3+(2D-3)\sqrt2]/(2D)$ for a particular $(2D+1)$-point set.
Zbarsky's cluster construction~\cite[Section~7]{Zbarsky14} gives the
larger limiting value $1.447317759362\ldots$ as the number of orthogonal clusters grows.
We give the finite correction and compare both bounds at fixed dimension,
including an explicit menu of variable simplex sizes.

We first state the vertex principle, then give the tree reduction,
joining lemma, and normalization transfer. The local estimates and a
direct coverage argument follow. The main induction is proved using the
finite base range $2\le m\le10$, whose independent proof then completes
the upper-bound argument. The terminal argument explicitly pays for the
connection to a one-point tail. Only after this do we present the two
lower-bound families and their dimension comparison.

\section{The vertex lemma}
\begin{lemma}[Vertex criterion]\label{lem:vertex}
Let $D\subseteq\mathbb R^m$ be a convex polytope with vertices
$P_1,\ldots,P_N$. Suppose $V:\mathbb R^m\to\mathbb R^r$ and
$g:\mathbb R^m\to\mathbb R$ are affine. If
$\norm{V(P_k)}\le g(P_k)$ for every vertex, then
$\norm{V(P)}\le g(P)$ for every $P\in D$.
\end{lemma}
\begin{proof}
Write $P=\sum_k\lambda_kP_k$ with $\lambda_k\ge0$ and
$\sum_k\lambda_k=1$. Affine maps preserve such combinations, so
\[
\norm{V(P)}=\norm{\sum_k\lambda_kV(P_k)}
\le\sum_k\lambda_k\norm{V(P_k)}
\le\sum_k\lambda_k g(P_k)=g(P).
\]
The first inequality is the triangle inequality.
\end{proof}
Every norm used below is the norm of a fixed affine vector in the retained
radii $d_4,\ldots,d_n$. Multiplication by a positive constant and shifts
inside a vector coordinate preserve the affine property. Thus the lemma
applies equally to vectors containing $d_i$ and to vectors containing
$d_i+$ a fixed constant.

The interval constraints and ordering inequalities defining the regions
below are linear, so their intersections and projections are convex
polytopes. For completeness, the ordered region
\[
a\le d_n\le d_{n-1}\le\cdots\le d_4\le b\qquad(a<b)
\]
has a particularly simple vertex description. Put $k=n-3$, and let
$P_j$, $0\le j\le k$, have its first $j$ coordinates equal to $b$ and
its remaining coordinates equal to $a$. For any point in the region set
\[
\lambda_0=\frac{b-d_4}{b-a},\qquad
\lambda_j=\frac{d_{j+3}-d_{j+4}}{b-a}\ (1\le j<k),\qquad
\lambda_k=\frac{d_n-a}{b-a}.
\]
Ordering makes these numbers nonnegative, and their sum telescopes to one.
In coordinate $\ell$, their weighted vertex sum is
\[
a+(b-a)\sum_{j=\ell}^k\lambda_j
=a+(d_{\ell+3}-a)=d_{\ell+3}.
\]
Thus $(d_4,\ldots,d_n)=\sum_{j=0}^k\lambda_jP_j$, proving directly that
these step profiles generate the whole region. They are its vertices:
the differences of consecutive profiles occupy distinct coordinate
directions, so the profiles are affinely independent. This justifies
the four-, five-, and seven-vertex lists used below.

\section{From rooted paths to trees and block induction}\label{sec:framework}
\subsection{Replacing child stars}
The following rooted-star replacement is the reduction used by
Khuller, Raghavachari, and Young~\cite[Section~4]{KRY96}; it is also
described by Chan~\cite[Section~2]{Chan04} and
Zbarsky~\cite[Section~2]{Zbarsky14}. We include its proof to make the
passage from paths to degree-three trees explicit.
\begin{proposition}[Rooted-path reduction]\label{prop:reduction}
If Corollary~\ref{cor:rooted} holds with a constant $C$, then
$\MSTthree(X)\le C\MST(X)$ for every finite Euclidean point set $X$.
\end{proposition}
\begin{proof}
Root a minimum spanning tree $T$ at an arbitrary vertex. At each vertex
$v$, replace the star joining $v$ to its children by a path starting at
$v$ and visiting all its child roots. The rooted-path hypothesis bounds
its length by $C$ times the sum of the original child-edge lengths.
Every edge of $T$ belongs to exactly one child star, so the sum of all
replacement lengths is at most $C\len(T)$.

The replacements form a tree. Inducting upward from the leaves, the
replacement at $v$ connects the already transformed child subtrees to
$v$, so the union is connected. If $v$ has $k_v$ children, its replacement
has $k_v$ edges, and $\sum_v k_v=|X|-1$. No two replacements use the same
edge: an edge joins either a parent to an original child or two original
siblings, and such a pair has a unique owner. A connected graph with
$|X|-1$ edges is a tree.

Each nonroot vertex has at most two incident edges in its parent's
replacement and at most one in its own replacement. Its degree is at
most three. The root has degree at most one. This gives the required
degree-three spanning tree.
\end{proof}

\subsection{The local extension problem}
Fix a block size $n\ge3$. The first $n$ labels form the block, and
$a_{n+1},a_{n+2}$ are its two incoming labels. An extension $Q_{s,j}$
starts at $a_s$, visits every block label once, and ends at $a_j$, where
$s\in\{n+1,n+2\}$ and $j\in\{1,2\}$. Its length contains $n$ edges.
A local estimate has the form
\begin{equation}\label{eq:localbudget}
\frac14\sum_{s=n+1}^{n+2}\sum_{j=1}^2\len(Q_{s,j})
\le C\sum_{i=1}^n d_i.
\end{equation}
The local paths start at incoming points, rather than at the origin.
The following is the joining step in Zbarsky's two-endpoint block
induction~\cite[Section~4]{Zbarsky14}, stated with the averages that we
will use throughout.
\begin{lemma}[Joining a block to a tail]\label{lem:join}
Suppose $m\ge n+2$. Assume that $T_{n+1},T_{n+2}$ start at the origin,
visit the actual tail labels $a_{n+1},\ldots,a_m$, end at the indicated
incoming labels, and satisfy
\[
\frac{\len(T_{n+1})+\len(T_{n+2})}{2}
\le C\sum_{i=n+1}^m d_i.
\]
If \eqref{eq:localbudget} holds, then the conclusion of
Theorem~\ref{thm:paths} holds for the full list.
\end{lemma}
\begin{proof}
For each $j$, concatenate $T_{n+1}$ with $Q_{n+1,j}$, or $T_{n+2}$
with $Q_{n+2,j}$, identifying the shared incoming label. Each candidate
starts at the origin, visits every actual label once, and ends at $a_j$.
Choose the shorter as $P_j$. Since a minimum is at most an average,
\[
\len(P_j)\le\frac12\sum_{s=n+1}^{n+2}
\bigl(\len(T_s)+\len(Q_{s,j})\bigr).
\]
Averaging over $j$ gives the sum of the tail-path average and the
four-extension average. Their budgets add to $C\sum_{i=1}^m d_i$.
\end{proof}

\subsection{Why the rooted-path constant cannot be below $3/2$}
The following example is recorded by Zbarsky~\cite[Section~3]{Zbarsky14}.
With root zero and points $1,-1$ on the real line, either rooted visiting
order has length three, whereas the radius sum is two. Thus any universal
rooted-path constant is at least $3/2$. This does not give a tree-ratio
lower bound of $3/2$: those same three collinear points already have an
MST of maximum degree two. The rooted-path problem is a sufficient
route to a degree-three bound, not an equivalent extremal problem.

\section{Normalization and its restoration cost}\label{sec:normalization}
The local estimates assume that the first three radii equal one. Scaling
justifies the value one once they are equal; a separate argument is
needed to make unequal radii equal. Radial contraction and restoration
are used in earlier analyses by Khuller, Raghavachari, and
Young~\cite[Sections~3--4]{KRY96}, Chan~\cite[Section~5]{Chan04}, and
Zbarsky~\cite[Section~4]{Zbarsky14}. The following lemma states and proves
the precise transfer estimate needed for our four-extension average.
\begin{lemma}[Normalization transfer]\label{lem:normalize}
Let $n\ge3$, $d_1\ge\cdots\ge d_{n+2}$, and $r=d_3>0$. Define
\[
\widetilde a_1=(r/d_1)a_1,\qquad
\widetilde a_2=(r/d_2)a_2,\qquad
\widetilde a_i=a_i\quad(i\ge3).
\]
If the normalized points $\widetilde a_i/r$ admit the four extensions
with a bound of the form \eqref{eq:localbudget} for some $C\ge3/2$,
then the original points admit that same bound.
\end{lemma}
\begin{proof}
Scale the normalized extensions back by $r$. Their average length is
at most $C\widetilde S$, where
\[
\widetilde S=2r+\sum_{i=3}^n d_i.
\]
Only labels 1 and 2 moved under contraction. Their restoration distances
are $\Delta_1=d_1-r$ and $\Delta_2=d_2-r$. Moving a path vertex by a
distance $h$ increases the path length by at most $h$ for an endpoint
and $2h$ for an internal vertex, by the triangle inequality on its
incident edges.

Across the four extensions, each moved label is the prescribed final
endpoint in two paths and is internal in the other two. It is never an
incoming endpoint because $n\ge3$. Its mean incidence is therefore
$(1+1+2+2)/4=3/2$. Restoring both labels, while preserving every visiting
order, gives a four-extension average of at most
\begin{align*}
C\widetilde S+\frac32(\Delta_1+\Delta_2)
&=C\sum_{i=1}^n d_i-\left(C-\frac32\right)(\Delta_1+\Delta_2)\\
&\le C\sum_{i=1}^n d_i.
\end{align*}
The starts remain fixed and the final labels remain $a_1,a_2$.
\end{proof}
Contraction may create coincident positions; the labeled-path convention
makes the argument valid without a perturbation. The selected block size
may depend on the normalized profile: the lemma transfers whichever
valid local estimate the cover supplies. If $d_3=0$, normalization is
unavailable, and the direct construction in Section~\ref{sec:base}
handles the configuration instead.

\section{A cover by eight ordered regions}\label{sec:cover}
Let $a_1,\ldots,a_{11}$ be labeled points in a Euclidean space, with
\[
d_i=\norm{a_i},\qquad d_1=d_2=d_3=1\ge d_4\ge\cdots\ge d_{11}\ge0.
\]
For a block size $n$, its vertices are $a_1,\ldots,a_n$ and its incoming
points are $a_{n+1},a_{n+2}$. For $s\in\{n+1,n+2\}$ and $j\in\{1,2\}$,
an extension $Q_{s,j}$ starts at $a_s$, visits every block vertex exactly
once, and ends at $a_j$. Set
\[
A_n=\frac14\sum_{s=n+1}^{n+2}\sum_{j=1}^2\len(Q_{s,j}),\qquad
S_n=3+\sum_{i=4}^n d_i,\qquad C=\frac{757}{500}=1.514.
\]
The sum is empty when $n=3$.

The eight regions below are all intersected with the standing ordering.
Any coordinate not restricted in the table has only its standing
normalization and ordering restrictions. The regions may overlap.
The subscript $i$ in $R_i$ is a region label; it is distinct from the
block size $n$. We list the regions in nondecreasing order of $n$.
There are two regions for $n=4$ and two for $n=6$, with different
radius conditions and different estimates.
\begin{center}\renewcommand{\arraystretch}{1.08}
\begin{tabular}{cclc}\toprule
Region&Block size $n$&Additional conditions&Vertex checks\\\midrule
$R_1$&3&$d_4\le1/4$&1\\
$R_2$&4&$1/4\le d_4\le1/2,\quad d_6\le1/4$&2\\
$R_3$&4&$d_4\ge1/2,\quad d_5\le3/8$&2\\
$R_4$&5&$d_4\ge1/2,\quad d_5\ge3/8,\quad d_6\le1/2$&4\\
$R_5$&6&$d_4\le1/2,\quad d_6\ge1/4$&4\\
$R_6$&6&$d_6\ge1/2,\quad d_7\le1/2$&4\\
$R_7$&7&$d_7\ge1/2,\quad d_9\le1/2$&5\\
$R_8$&9&$d_9\ge1/2$&7\\
\bottomrule
\end{tabular}\end{center}
\begin{theorem}\label{thm:local}
Every normalized ordered radius profile belongs to at least one of
$R_1,\ldots,R_8$. On each $R_i$, with its assigned block size $n$,
the four extensions can be chosen so that $A_n\le C S_n$.
\end{theorem}
The proof uses the rational path mixtures and positive witness vectors given below. In particular, the block sizes $n=8$
and $n=10$ are unnecessary for this cover. We do not claim that eight is
the smallest possible number of regions.

\subsection{A direct proof of coverage}\label{sec:coverage}
The following decision procedure selects a valid region for every
profile. Its branch order need not agree with the numerical order of
the region labels. At a threshold equality, either applicable branch
may be used.
First suppose $d_4\le1/2$.
If $d_4\le1/4$, use $R_1$. Otherwise $1/4\le d_4\le1/2$.
If $d_6\le1/4$, use $R_2$; if $d_6\ge1/4$, use $R_5$.
This covers every profile with $d_4\le1/2$.

Now suppose $d_4\ge1/2$.
If $d_5\le3/8$, use $R_3$. Otherwise $d_5\ge3/8$.
If $d_6\le1/2$, use $R_4$. Otherwise $d_6\ge1/2$.
If $d_7\le1/2$, use $R_6$. Otherwise $d_7\ge1/2$.
Finally, use $R_7$ when $d_9\le1/2$ and $R_8$ when $d_9\ge1/2$.
These alternatives exhaust all profiles. Since every region is closed,
all threshold equalities are included. This is a proof on the continuous
domain, rather than a verification on sampled profiles.

\section{Path averages and the weighted-distance estimate}\label{sec:weighted}
\subsection{Path mixtures}
Averaging admissible visiting orders and applying weighted-distance
witnesses is the method used by Zbarsky~\cite[Sections~3--4]{Zbarsky14}.
The template probabilities and witness choices below specify the local
estimates used in this paper.
For $n\in\{4,5,6,7,9\}$ and either incoming point $b$, use the templates
\[
\begin{array}{ccl}
\mathrm I&:&b,a_1,a_3,a_4,\ldots,a_n,a_2,\\
\mathrm {II}&:&b,a_3,a_1,a_4,\ldots,a_n,a_2,\\
\mathrm {III}&:&b,a_3,a_4,\ldots,a_n,a_1,a_2.
\end{array}
\]
Choose a template with the following probabilities, then independently
interchange labels 1 and 2 uniformly and permute labels $3,\ldots,n$
uniformly. The permutations affect the visiting order, not the points.
\begin{center}\begin{tabular}{cccc}\toprule
$n$&$p_n$ (I)&$q_n$ (II)&$r_n$ (III)\\\midrule
4&$2/5$&$3/10$&$3/10$\\
5&$1/3$&$5/12$&$1/4$\\
6&$7/25$&$13/25$&$1/5$\\
7&$7/30$&$3/5$&$1/6$\\
9&$7/40$&$7/10$&$1/8$\\\bottomrule\end{tabular}\end{center}
Each family contains $2(n-2)!$ paths and gives equal total weight to
paths ending at $a_1$ and at $a_2$. Thus the full mixture assigns weight
$1/2$ to each endpoint class. Let $E_b$ be its mean length.
Choosing a shortest listed path for each endpoint gives
\begin{equation}\label{eq:select}
\frac{\len(Q_{s,1})+\len(Q_{s,2})}{2}\le E_{a_s},\qquad
A_n\le\frac{E_{a_{n+1}}+E_{a_{n+2}}}{2}.
\end{equation}
The $n=3$ mixture is given in Section~\ref{sec:r1}.

\subsection{Edge weights and the witness inequality}
The edge families are terminal--terminal, terminal--middle,
middle--middle, incoming--terminal, and incoming--middle, where
$a_1,a_2$ are terminals and $a_3,\ldots,a_n$ are middle vertices.
Their respective weights per edge are
\[
r_n,\quad \frac{2p_n+3q_n+r_n}{2(n-2)},\quad
\frac{(n-3)(p_n+r_n)+(n-4)q_n}{\binom{n-2}{2}},\quad
\frac{p_n}{2},\quad\frac{q_n+r_n}{n-2}.
\]
These follow by counting edge occurrences in the three templates and
then using symmetry within a family. Explicitly:
\begin{center}\small\begin{tabular}{cccccc}\toprule
$n$&Terminal--terminal&Terminal--middle&Middle--middle&Incoming--terminal&Incoming--middle\\\midrule
4&$3/10$&$1/2$&$7/10$&$1/5$&$3/10$\\
5&$1/4$&$13/36$&$19/36$&$1/6$&$2/9$\\
6&$1/5$&$29/100$&$31/75$&$7/50$&$9/50$\\
7&$1/6$&$73/300$&$17/50$&$7/60$&$23/150$\\
9&$1/8$&$103/560$&$53/210$&$7/80$&$33/280$\\\bottomrule
\end{tabular}\end{center}
\subsection{Young's weighted-distance inequality}
The following is the upper-bound inequality in Young's dual formulation
\cite[Lemma~1]{Young94}, also used in Zbarsky's analysis
\cite[Section~3]{Zbarsky14}. We require only this inequality, and no
assumption that the witness is optimal.
\begin{theorem}[Young's weighted-distance inequality]\label{thm:young}
Let $y_1,\ldots,y_k$ be vectors in a Euclidean space, let
$w_{ij}\ge0$ for $1\le i<j\le k$, and let $\ell_i\ge\norm{y_i}$.
For every choice of positive numbers $x_1,\ldots,x_k$,
\begin{equation}\label{eq:young}
\sum_{i<j}w_{ij}\norm{y_i-y_j}
\le
\left[
\left(\sum_{i<j}\frac{w_{ij}^2}{x_ix_j}\right)
\left(\sum_i x_i\right)
\left(\sum_i x_i\ell_i^2\right)
\right]^{1/2}.
\end{equation}
In particular, one may take $\ell_i=\norm{y_i}$.
\end{theorem}
\begin{proof}
Apply the Cauchy--Schwarz inequality to the lists indexed by unordered
pairs,
\[
\frac{w_{ij}}{\sqrt{x_ix_j}}
\quad\hbox{and}\quad
\sqrt{x_ix_j}\,\norm{y_i-y_j}.
\]
This gives
\begin{equation}\label{eq:young-cs}
\left(\sum_{i<j}w_{ij}\norm{y_i-y_j}\right)^2
\le
\left(\sum_{i<j}\frac{w_{ij}^2}{x_ix_j}\right)
\left(\sum_{i<j}x_ix_j\norm{y_i-y_j}^2\right).
\end{equation}
The second factor has the exact expansion
\begin{align*}
\sum_{i<j}x_ix_j\norm{y_i-y_j}^2
&=\frac12\sum_{i,j}x_ix_j
  \bigl(\norm{y_i}^2+\norm{y_j}^2-2y_i\cdot y_j\bigr)\\
&=\left(\sum_i x_i\right)\left(\sum_i x_i\norm{y_i}^2\right)
  -\norm{\sum_i x_i y_i}^2\\
&\le\left(\sum_i x_i\right)\left(\sum_i x_i\ell_i^2\right).
\end{align*}
Here the squared norm is nonnegative and $\norm{y_i}\le\ell_i$.
Substitute this bound into \eqref{eq:young-cs} and take nonnegative
square roots to obtain \eqref{eq:young}.
\end{proof}

\subsection{Applying Young's inequality to the two incoming points}
With $(y_1,\ldots,y_{n+1})=(a_1,\ldots,a_n,b)$, denote the edge weights
from the preceding table by $w_{ij}$. Then
$E_b=\sum_{i<j}w_{ij}\norm{y_i-y_j}$. For a positive witness set
\[
Q_x=\sum_{i<j}\frac{w_{ij}^2}{x_ix_j},\qquad
X_x=\sum_i x_i,\qquad K_x=Q_xX_x.
\]
Theorem~\ref{thm:young}, with $\ell_i=\norm{y_i}$, gives
\[
E_b\le\sqrt{K_x\sum_i x_i\norm{y_i}^2}.
\]
All witnesses below use the same entries for the two incoming points.
Define the auxiliary quantity
\[
\eta_n=\frac{d_{n+1}^2+d_{n+2}^2}{2}.
\]
Averaging the two incoming estimates in \eqref{eq:select} and using
$(\sqrt R+\sqrt S)/2\le\sqrt{(R+S)/2}$ proves
\begin{equation}\label{eq:master}
A_n\le\sqrt{K_x\left(\sum_{i=1}^n x_i d_i^2+x_{n+1}\eta_n\right)}.
\end{equation}
This step bounds the two incoming means together; it does not replace
either incoming squared radius individually by their average.

\subsection{How the checks are made}
The following sections specify $V_i$ with $A_n\le\norm{V_i}$ on $R_i$.
The retained coordinates are always denoted by their original names
$d_4,\ldots,d_n$. Each $V_i$ is affine. The right-hand target $C S_n$
is affine as well, so Lemma~\ref{lem:vertex} reduces each comparison to
the listed vertices.

The left-hand numerical column in each table is rounded \emph{upwards}
to six decimal places; it is a rational upper bound for $\norm{V_i}$.
The right-hand numerical column is the exact value of $C S_n$.
For each row, squaring the displayed rational upper bound and comparing
it with the squared norm given by the formula verifies the rounding
direction using integer arithmetic. Each upper bound is strictly smaller
than its right-hand target. Thus the printed tables give rigorous finite
comparisons, without an external numerical certificate or a square-root
approximation used to infer a sign.

\section{Verification of the eight local estimates}\label{sec:cases}

\subsection{Region $R_1$: $n=3$}\label{sec:r1}
Here $d_4\le1/4$, hence $d_5\le1/4$. For either incoming point, use
\[
\begin{array}{c|cccc}
\text{Order}&(b,a_1,a_3,a_2)&(b,a_2,a_3,a_1)&(b,a_3,a_1,a_2)&(b,a_3,a_2,a_1)\\\hline
\text{Weight}&9/34&9/34&4/17&4/17.
\end{array}
\]
Each endpoint receives weight $1/2$. The weights of
$a_1a_2$, $a_1a_3,a_2a_3$, $ba_1,ba_2$, and $ba_3$ are respectively
$8/17$, $13/17$, $9/34$, and $8/17$.

\paragraph{Witness.} In the order $(a_1,\ldots,a_{3},b)$, use
\[x=(1,1,7/5,14/15),\qquad K_x=\frac{1013597}{169932}.\]
The radius factor is $17/5+(14/15)\eta_3$. Since $\eta_3\le1/16$,
\[
V_1=\left(\sqrt{\frac{1013597}{169932}\cdot\frac{83}{24}}\right),
\qquad A_3\le\norm{V_1}.
\]
There are no variable block radii. The retained domain is a singleton,
denoted by the empty coordinate tuple $()$.

\begin{center}\small\renewcommand{\arraystretch}{1.08}
\begin{tabular}{crr}\toprule
Corner $()$ & $\norm{V_i}\le$ & $C S_n$ \\\midrule
$()$ & 4.541806 & 4.542000\\
\bottomrule\end{tabular}\end{center}

\subsection{Region $R_2$: $n=4$}
Here $1/4\le d_4\le1/2$ and $d_6\le1/4$.
The first incoming radius $d_5$ is allowed to be as large as $d_4$.

\paragraph{Witness.} In the order $(a_1,\ldots,a_{4},b)$, use
\[x=(1,1,7/5,12/5,1),\qquad K_x=\frac{70193}{10500}.\]
Using $d_5\le d_4$ and $d_6\le1/4$, we have
\[
\eta_4\le\frac12d_4^2+\frac1{32}.
\]
Therefore \eqref{eq:master} gives $A_4\le\norm{V_2}$, where
\[
V_2=\sqrt{\frac{70193}{10500}}
\left(\sqrt{\frac{549}{160}},\sqrt{\frac{29}{10}}d_4\right).
\]
The retained domain is the interval $[1/4,1/2]$.

\begin{center}\small\renewcommand{\arraystretch}{1.08}
\begin{tabular}{crr}\toprule
Corner $d_4$ & $\norm{V_i}\le$ & $C S_n$ \\\midrule
$(1/4)$ & 4.914238 & 4.920500\\
$(1/2)$ & 5.271123 & 5.299000\\
\bottomrule\end{tabular}\end{center}

\subsection{Region $R_3$: $n=4$}
Here $d_4\ge1/2$ and $d_5\le3/8$.

\paragraph{Witness.} In the order $(a_1,\ldots,a_{4},b)$, use
\[x=(1,1,7/5,8/5,6/7),\qquad K_x=\frac{334027}{47040}.\]
Both incoming radii $d_5,d_6$ are at most $3/8$, so $\eta_4\le9/64$.
Thus $A_4\le\norm{V_3}$ with
\[
V_3=\sqrt{\frac{334027}{47040}}
\left(\sqrt{\frac{17}{5}+\frac67\cdot\frac9{64}},\sqrt{\frac85}d_4\right).
\]
The retained domain is $1/2\le d_4\le1$.

\begin{center}\small\renewcommand{\arraystretch}{1.08}
\begin{tabular}{crr}\toprule
Corner $d_4$ & $\norm{V_i}\le$ & $C S_n$ \\\midrule
$(1/2)$ & 5.276305 & 5.299000\\
$(1)$ & 6.029966 & 6.056000\\
\bottomrule\end{tabular}\end{center}

\newpage
\subsection{Region $R_4$: $n=5$}
Here $d_4\ge1/2$, $d_5\ge3/8$, and $d_6\le1/2$.

\paragraph{Witness.} In the order $(a_1,\ldots,a_{5},b)$, use
\[x=(1,1,11/8,8/5,9/5,1),\qquad K_x=\frac{333812783}{41057280}.\]
The master estimate is
\[
A_5\le\sqrt{\frac{333812783}{41057280}
\left(\frac{27}{8}+\frac85d_4^2+\frac95d_5^2+\eta_5\right)}.
\]
Ordering and $d_6\le1/2$ imply
\[
\eta_5\le\min\{d_5^2,1/4\}.
\]
We use one quadratic majorant for this minimum, so no subdivision of the
region is needed. Choose a quadratic agreeing with $d_5^2$ at $d_5=3/8$
and $d_5=1/2$, and with $1/4$ at $d_5=1$. Solving these three linear
conditions on its coefficients gives
\[
q(d_5)=-\frac9{20}+\frac{21}{10}d_5-\frac75d_5^2.
\]
For $3/8\le d_5\le1/2$,
\[
q(d_5)-d_5^2=-\frac{12}{5}(d_5-3/8)(d_5-1/2)\ge0.
\]
For $1/2\le d_5\le1$,
\[
q(d_5)-1/4=-\frac75(d_5-1/2)(d_5-1)\ge0.
\]
These two sign calculations prove $\eta_5\le q(d_5)$ throughout $R_4$.
Completing the square gives the single formula
\[
\frac{27}{8}+\frac85d_4^2+\frac95d_5^2+q(d_5)
=\frac{27}{160}+\frac85d_4^2+\frac25(d_5+21/8)^2.
\]
Hence $A_5\le\norm{V_4}$, where
\[
V_4=\sqrt{\frac{333812783}{41057280}}
\left(\sqrt{\frac{27}{160}},\sqrt{\frac85}d_4,
\sqrt{\frac25}(d_5+21/8)\right).
\]
This vector is affine in $(d_4,d_5)$. The retained region
\[
1/2\le d_4\le1,\qquad 3/8\le d_5\le d_4
\]
is a quadrilateral with the four vertices below. The quadratic majorant
agrees with the minimum at $d_5=3/8,1/2,1$, which are precisely the
$d_5$-coordinates of these vertices.

\begin{center}\small\renewcommand{\arraystretch}{1.08}
\begin{tabular}{crr}\toprule
Corner $(d_4,d_5)$ & $\norm{V_i}\le$ & $C S_n$ \\\midrule
$(1/2,3/8)$ & 5.821828 & 5.866750\\
$(1/2,1/2)$ & 6.031884 & 6.056000\\
$(1,3/8)$ & 6.606828 & 6.623750\\
$(1,1)$ & 7.557525 & 7.570000\\
\bottomrule\end{tabular}\end{center}

\subsection{Region $R_5$: $n=6$}
The conditions $d_4\le1/2$ and $d_6\ge1/4$ give
\[
1/4\le d_6\le d_5\le d_4\le1/2.
\]

\paragraph{Witness.} In the order $(a_1,\ldots,a_{6},b)$, use
\[x=(1,1,11/8,5/2,5/2,5/2,5/4),\qquad K_x=\frac{135359911}{16500000}.\]
Since $d_8\le d_7\le d_6$, we have $\eta_6\le d_6^2$.
Consequently $A_6\le\norm{V_5}$ with
\[
V_5=\sqrt{\frac{135359911}{16500000}}
\left(\sqrt{\frac{27}{8}},\sqrt{\frac52}d_4,
\sqrt{\frac52}d_5,\sqrt{\frac{15}{4}}d_6\right).
\]
The retained ordered tetrahedron has the following four vertices.

\begin{center}\small\renewcommand{\arraystretch}{1.08}
\begin{tabular}{crr}\toprule
Corner $(d_4,d_5,d_6)$ & $\norm{V_i}\le$ & $C S_n$ \\\midrule
$(1/4,1/4,1/4)$ & 5.672180 & 5.677500\\
$(1/2,1/4,1/4)$ & 6.001589 & 6.056000\\
$(1/2,1/2,1/4)$ & 6.313836 & 6.434500\\
$(1/2,1/2,1/2)$ & 6.755198 & 6.813000\\
\bottomrule\end{tabular}\end{center}

\subsection{Region $R_6$: $n=6$}
Here $d_6\ge1/2$ and $d_7\le1/2$, so the retained domain is
$1/2\le d_6\le d_5\le d_4\le1$.

\paragraph{Witness.} In the order $(a_1,\ldots,a_{6},b)$, use
\[x=(1,1,11/8,8/5,8/5,8/5,9/10),\qquad K_x=\frac{22828003}{2400000}.\]
Both incoming radii $d_7,d_8$ are at most $1/2$.
Thus $\eta_6\le1/4$, and $A_6\le\norm{V_6}$ with
\[
V_6=\sqrt{\frac{22828003}{2400000}}
\left(\sqrt{\frac{18}{5}},\sqrt{\frac85}d_4,
\sqrt{\frac85}d_5,\sqrt{\frac85}d_6\right).
\]

\begin{center}\small\renewcommand{\arraystretch}{1.08}
\begin{tabular}{crr}\toprule
Corner $(d_4,d_5,d_6)$ & $\norm{V_i}\le$ & $C S_n$ \\\midrule
$(1/2,1/2,1/2)$ & 6.756923 & 6.813000\\
$(1,1/2,1/2)$ & 7.554470 & 7.570000\\
$(1,1,1/2)$ & 8.275507 & 8.327000\\
$(1,1,1)$ & 8.938569 & 9.084000\\
\bottomrule\end{tabular}\end{center}

\subsection{Region $R_7$: $n=7$}
Here $d_7\ge1/2$ and $d_9\le1/2$. Importantly, the first incoming
radius $d_8$ need not be at most $1/2$: it can be as large as $d_7$.

\paragraph{Witness.} In the order $(a_1,\ldots,a_{7},b)$, use
\[x=(1,1,11/8,5/3,5/3,5/3,5/3,1),\qquad K_x=\frac{207766943}{19800000}.\]
Since $d_8\le d_7$ and $d_9\le1/2$,
\[
\eta_7\le\frac12d_7^2+\frac18.
\]
Consequently $A_7\le\norm{V_7}$ with
\[
V_7=\sqrt{\frac{207766943}{19800000}}
\left(\sqrt{\frac72},\sqrt{\frac53}d_4,
\sqrt{\frac53}d_5,\sqrt{\frac53}d_6,\sqrt{\frac{13}{6}}d_7\right).
\]
The retained domain is $1/2\le d_7\le\cdots\le d_4\le1$.
The five step vertices below suffice. Allowing this larger first incoming
radius is what makes a separate $n=8$ rule unnecessary for coverage.

\begin{center}\small\renewcommand{\arraystretch}{1.08}
\begin{tabular}{crr}\toprule
Corner $(d_4,d_5,d_6,d_7)$ & $\norm{V_i}\le$ & $C S_n$ \\\midrule
$(1/2,1/2,1/2,1/2)$ & 7.451641 & 7.570000\\
$(1,1/2,1/2,1/2)$ & 8.285140 & 8.327000\\
$(1,1,1/2,1/2)$ & 9.042132 & 9.084000\\
$(1,1,1,1/2)$ & 9.740470 & 9.841000\\
$(1,1,1,1)$ & 10.579619 & 10.598000\\
\bottomrule\end{tabular}\end{center}

\subsection{Region $R_8$: $n=9$}
Here $d_9\ge1/2$. The only incoming restrictions needed are the
standing ordering inequalities $d_{11}\le d_{10}\le d_9$.

\paragraph{Witness.} In the order $(a_1,\ldots,a_{9},b)$, use
\[x=(1,1,11/8,8/5,8/5,8/5,8/5,8/5,8/5,9/10),\qquad K_x=\frac{2186251079}{165580800}.\]
Thus $\eta_9\le d_9^2$, and $A_9\le\norm{V_8}$, where
\[
V_8=\sqrt{\frac{2186251079}{165580800}}
\left(\sqrt{\frac{27}{8}},\sqrt{\frac85}d_4,\ldots,
\sqrt{\frac85}d_8,\sqrt{\frac52}d_9\right).
\]
The retained domain is $1/2\le d_9\le\cdots\le d_4\le1$.
It has seven step vertices. This region includes every normalized ordered
profile with $d_{10}\ge1/2$, so a separate $n=10$ rule is unnecessary.

\begin{center}\small\renewcommand{\arraystretch}{1.08}
\begin{tabular}{crr}\toprule
Corner $(d_4,d_5,d_6,d_7,d_8,d_9)$ & $\norm{V_i}\le$ & $C S_n$ \\\midrule
$(1/2,1/2,1/2,1/2,1/2,1/2)$ & 8.900629 & 9.084000\\
$(1,1/2,1/2,1/2,1/2,1/2)$ & 9.750150 & 9.841000\\
$(1,1,1/2,1/2,1/2,1/2)$ & 10.531366 & 10.598000\\
$(1,1,1,1/2,1/2,1/2)$ & 11.258503 & 11.355000\\
$(1,1,1,1,1/2,1/2)$ & 11.941446 & 12.112000\\
$(1,1,1,1,1,1/2)$ & 12.587389 & 12.869000\\
$(1,1,1,1,1,1)$ & 13.535102 & 13.626000\\
\bottomrule\end{tabular}\end{center}

\section{The global induction and the spanning-tree bound}\label{sec:global}
The local work proves Theorem~\ref{thm:local}: coverage was established
in Section~\ref{sec:coverage}, and the path mixtures and the 29 vertex
checks establish the four-extension bounds on its eight regions.

\begin{proof}[Proof of Theorem~\ref{thm:paths}, using the base cases below]
Section~\ref{sec:base} proves the result for $2\le m\le10$ using only
the local estimates, normalization, and direct terminal constructions.
Assume now $m\ge11$ and that the theorem holds for every smaller number
of points at least two.

If $d_3=0$, apply Lemma~\ref{lem:zeros} below. Otherwise contract the
first two radii to $d_3$ and scale by $1/d_3$. The resulting profile
satisfies the normalization for the eight-region cover. Some region
supplies a block size $n\in\{3,4,5,6,7,9\}$ and the four normalized
extensions. Lemma~\ref{lem:normalize} restores the original radii and
establishes \eqref{eq:localbudget}.

The actual tail has $m-n\ge11-9=2$ labels and fewer than $m$ labels.
The strong induction hypothesis gives its two rooted paths, ending at
the two largest-radius tail labels $a_{n+1},a_{n+2}$.
Lemma~\ref{lem:join} combines them with the four local extensions.
This proves the inductive step.
\end{proof}
Corollary~\ref{cor:rooted} follows by translation and taking the shorter
of the two paths, with the zero- and one-point cases handled directly.
Proposition~\ref{prop:reduction} then proves Theorem~\ref{thm:tree}.

\subsection{Constructivity and the extent of the result}
As in Zbarsky's bounded-block reconstruction~\cite[Section~6]{Zbarsky14},
the argument is constructive. The cover is selected by a short sequence
of radius comparisons. Every local mixture contains finitely many paths
of bounded size, at most $3\cdot2\cdot7!=30{,}240$ for a nine-point block.
For each incoming label and prescribed endpoint, one may enumerate the
listed family and choose a shortest path. Positive mixture weights are
used to justify its length, not to randomize the construction. Each
recursive step removes at least three labels, so the recursion has a
linear number of block steps. Together with sorting the radii and the
star replacement, this gives a deterministic polynomial-time
construction once an MST is supplied.

The root-to-path constant proved here lies between $3/2$ and $757/500$.
We do not establish that the latter value, the witnesses, or the number
of regions is optimal. Nor does the lower limit $3/2$ for rooted paths
assert the same lower limit for degree-three spanning trees. The new
upper bound concerns total Euclidean edge length in arbitrary dimension;
it is separate from bottleneck objectives and bounds that use a fixed
planar embedding.

\section{Terminal constructions and all base cases}\label{sec:base}
This section supplies the base range invoked in Section~\ref{sec:global}.
Its proof uses the already established local estimates, coverage, and
normalization, and no instance of the global theorem outside the base
range. We count actual labels, including actual points at the origin;
only artificial padding labels may be removed.

\subsection{Two points and zero-radius labels}
For two labels use
\[
P_1:0,a_2,a_1,\qquad P_2:0,a_1,a_2.
\]
Their average length is
\begin{equation}\label{eq:two}
\frac{d_1+d_2}{2}+\norm{a_1-a_2}
\le\frac32(d_1+d_2)\le C(d_1+d_2).
\end{equation}
\begin{lemma}\label{lem:zeros}
For every $m\ge2$, if $d_3=\cdots=d_m=0$, the two-endpoint conclusion
holds with constant $3/2$.
\end{lemma}
\begin{proof}
For $m=2$ use \eqref{eq:two}. Otherwise take the orders
\[
0,a_3,\ldots,a_m,a_2,a_1,
\qquad 0,a_3,\ldots,a_m,a_1,a_2.
\]
The initial edges have length zero, and the remaining lengths are those
in \eqref{eq:two}. This includes configurations with fewer than two
positive radii.
\end{proof}

\subsection{Padding and the selected block}
Zero-padding is also used by Zbarsky~\cite[Section~4]{Zbarsky14}.
For selecting a region, append labeled copies of the origin until the
list has eleven labels if necessary. Order all actual labels before
padding labels at radius zero. If $d_3>0$, contract and scale as in
Lemma~\ref{lem:normalize}. The eight-region cover then applies.

Every assigned block has a positive smallest radius: $d_3=1$ on $R_1$,
$d_4\ge1/4$ on $R_2$, $d_4\ge1/2$ on $R_3$,
$d_5\ge3/8$ on $R_4$, $d_6\ge1/4$ on $R_5$,
$d_6\ge1/2$ on $R_6$, $d_7\ge1/2$ on $R_7$, and
$d_9\ge1/2$ on $R_8$.
Consequently the selected block contains only actual positive-radius
labels, and $n\le m$. Its incoming labels may include padding when
$m-n$ equals zero or one. Apply the local estimate and normalization
transfer before carrying out either of the following constructions.

\subsection{A block leaving no actual tail point}
Suppose $m=n$. Let $z_1,z_2$ be the artificial incoming labels at the
origin. The restored local estimate gives
\[
\frac14\sum_{s=1}^2\sum_{j=1}^2\len(Q_{z_s,j})
\le C\sum_{i=1}^m d_i.
\]
For each $j$, take the shorter of the two extensions ending at $a_j$
and identify its initial artificial label with the designated root.
The resulting path visits every actual label once. The average of the
two selected lengths is at most the displayed four-extension average.
There is no additional tail cost.

\subsection{A block leaving one actual tail point}
Suppose $m=n+1$. Write $b=a_m$ and $r=\norm b$, and let $z$ be the
artificial incoming label at the origin. The four local extensions satisfy
\begin{equation}\label{eq:onetail-local}
\frac14\sum_{j=1}^2
\bigl(\len(Q_{b,j})+\len(Q_{z,j})\bigr)
\le C\sum_{i=1}^n d_i.
\end{equation}
For an endpoint $a_j$, first consider the root edge $0,b$ followed by
$Q_{b,j}$. It has length $r+\len(Q_{b,j})$.

A second candidate follows $0,b,z$ and then $Q_{z,j}$, with length
$2r+\len(Q_{z,j})$ before deleting $z$.
If the first block label after $z$ is $a_\ell$, replace the segment
sequence $b,z,a_\ell$ by $b,a_\ell$. The triangle inequality shows that
this shortcut cannot increase length. It removes only the artificial
label, so the candidate still visits every actual label once and ends
at $a_j$.

Choose the shorter candidate as $P_j$. Then
\[
\len(P_j)\le\frac{3r+\len(Q_{b,j})+\len(Q_{z,j})}{2}.
\]
Averaging over $j$ and applying \eqref{eq:onetail-local} gives
\begin{align*}
\frac{\len(P_1)+\len(P_2)}2
&\le\frac32r+\frac14\sum_{j=1}^2
\bigl(\len(Q_{b,j})+\len(Q_{z,j})\bigr)\\
&\le\frac32r+C\sum_{i=1}^n d_i
\le C\sum_{i=1}^m d_i.
\end{align*}
The final inequality uses $C\ge3/2$. In particular, introducing an
artificial incoming origin has not treated the one-point tail as free;
its full possible connection cost has been included.

\subsection{The first nontrivial sizes}
For $m=3$, the case $d_3=0$ is already settled. Otherwise normalization
makes the first three radii one and the padded $d_4,d_5$ zero.
Region $R_1$ supplies a three-point block with no actual tail. Transfer
the four extensions back by Lemma~\ref{lem:normalize}, then apply the
no-tail construction.

For $m=4$ with $d_3>0$, consider the normalized $d_4\in[0,1]$.
If $d_4\le1/4$, use $R_1$ with $n=3$ and the one-tail construction.
If $1/4\le d_4\le1/2$, use $R_2$ with $n=4$, since padded $d_6=0$.
If $d_4\ge1/2$, use $R_3$ with $n=4$, since padded $d_5=0$.
The latter two choices have no actual tail. The normalization lemma
transfers every estimate back to the original radii.

\subsection{The complete finite induction}
\begin{proposition}\label{prop:base}
Theorem~\ref{thm:paths} holds for every $m\in\{2,\ldots,10\}$.
\end{proposition}
\begin{proof}
Use strong induction within this finite range. The case $m=2$ is
\eqref{eq:two}. Fix $3\le m\le10$ and assume all smaller cases at least
two. If $d_3=0$, use Lemma~\ref{lem:zeros}. Otherwise pad, normalize,
and select a region. Its block size satisfies
$n\in\{3,4,5,6,7,9\}$ and $n\le m$, as proved above.
The local estimate transfers to the original four extensions by
Lemma~\ref{lem:normalize}.

Set $k=m-n$. If $k=0$, use the no-tail construction. If $k=1$, use the
one-tail construction. If $k\ge2$, then $k<m$, so the finite induction
hypothesis supplies the two actual tail paths, and Lemma~\ref{lem:join}
completes the configuration. These alternatives exhaust every $k\ge0$.
\end{proof}
For clarity, the permitted block sizes and corresponding tail sizes in
this finite induction are listed below. A listed size may be excluded
by a particular profile's region conditions; the table records all
sizes that can be selected at that value of $m$.
\begin{center}\small\renewcommand{\arraystretch}{1.08}
\begin{tabular}{ccl}\toprule
Actual labels $m$&Possible block sizes $n$&Tail sizes $m-n$\\\midrule
2&Direct construction&No block\\
3&$3$&$0$\\
4&$3,4$&$0,1$\\
5&$3,4,5$&$0,1,2$\\
6&$3,4,5,6$&$0,1,2,3$\\
7&$3,4,5,6,7$&$0,1,2,3,4$\\
8&$3,4,5,6,7$&$1,2,3,4,5$\\
9&$3,4,5,6,7,9$&$0,2,3,4,5,6$\\
10&$3,4,5,6,7,9$&$1,3,4,5,6,7$\\\bottomrule
\end{tabular}\end{center}
The largest terminal one-tail case is now $m=10,n=9$. The two candidates
above pay at most $(3/2)d_{10}$ in average connection cost in addition
to the nine-point block budget. All other actual tails of size at least
two use a strictly smaller, already proved case. This completes the
base range without using the global induction step.

\section{Lower bounds and comparison by dimension}\label{sec:lower}
The upper theorem concerns all finite configurations. A lower-bound
construction gives a particular point set for which every degree-three
tree is long relative to an unrestricted spanning tree. Following the ratio notation used by Chan~\cite{Chan04} and
Zbarsky~\cite[Section~2]{Zbarsky14}, for $D\ge2$ let
\[
\tau_3^D=\sup_{X\subset\mathbb R^D}
\frac{\MSTthree(X)}{\MST(X)},\qquad
\tau_3^\infty=\sup_{D\ge1}\tau_3^D,
\]
where the first supremum is over finite sets with at least two points.
Embedding by zero coordinates shows that $\tau_3^D$ is nondecreasing in $D$.

We give the two constructions considered here: a cross-polytope chain
within each fixed dimension, and orthogonal simplex clusters. The
seven-leaf instance of the second construction is Zbarsky's
construction~\cite[Section 7]{Zbarsky14}. The variable-size version below
uses the same idea and charging argument. Its purpose is to make the
finite-dimension comparison explicit; it does not improve Zbarsky's
limiting lower bound.

In proving a ratio lower bound, an upper bound on the unrestricted MST
length is sufficient. Indeed, if every degree-three tree has length at
least $U>0$ and some ordinary spanning tree has length $B>0$, then
\[
\frac{\MSTthree(X)}{\MST(X)}\ge\frac{U}{\MST(X)}\ge\frac{U}{B}.
\]
The first construction will determine both optimum lengths exactly in
its single-cell instance; the second needs only the displayed spanning
tree as a denominator bound.

\subsection{A cross-polytope chain in fixed dimension}
Fix integers $D\ge2$ and $N\ge1$, and let $e_1,\ldots,e_D$ be the
coordinate unit vectors. Take the centers, axial points, and transverse
points
\begin{align*}
p_i&=2i e_1 &&(1\le i\le N),\\
b_i&=(2i+1)e_1 &&(0\le i\le N),\\
p_i+e_j,\quad p_i-e_j&&&(1\le i\le N,\ 2\le j\le D).
\end{align*}
There are $2DN+1$ distinct integer-coordinate points. Join each $p_i$ to
$b_{i-1},b_i$ and to its $2(D-1)$ transverse points. These $2DN$ unit
edges form a tree: the axial points and centers make the alternating path
\[
b_0,p_1,b_1,p_2,\ldots,p_N,b_N,
\]
and every transverse point is a leaf. Every pair of distinct integer
points is at distance at least one. Since any spanning tree has $2DN$
edges, this proves
\begin{equation}\label{eq:chainmst}
\MST(X)=2DN.
\end{equation}

Every possible unit edge is incident with a center. To check this,
integer vectors at distance one must differ by $\pm1$ in exactly one
coordinate. Distinct axial points differ by at least two in their first
coordinate. An axial point and a transverse point differ in both the
first and a transverse coordinate. Two transverse points at different
centers differ by at least two in the first coordinate; at the same
center, their difference has either one coordinate of magnitude two or
two coordinates of magnitude one. Centers and axial points have no
nonzero transverse coordinate. These observations leave only the unit
edges incident with a center.

In a degree-three tree the $N$ centers have total degree at most $3N$,
so at most $3N$ edges have length one. Every remaining edge has squared
integer length at least two, hence length at least $\sqrt2$.
If $u\le3N$ is the number of unit edges, then
\[
\len(T)\ge u+(2DN-u)\sqrt2
\ge3N+(2D-3)N\sqrt2,
\]
because $1-\sqrt2<0$. Dividing by \eqref{eq:chainmst} proves the bound
\begin{proposition}\label{prop:chain}
For every $D\ge2$,
\[
\tau_3^D\ge L_D:=\frac{3+(2D-3)\sqrt2}{2D}
=\sqrt2-\frac{3(\sqrt2-1)}{2D}.
\]
\end{proposition}
For $D=2$, the bound $L_2=(3+\sqrt2)/4$ is the classical
square-plus-center lower bound of Khuller, Raghavachari, and
Young~\cite[Section~5]{KRY96}; see also Chan~\cite[Section~6]{Chan04}.
The coordinate construction above extends this example to all $D\ge2$.

The dimension stays equal to $D$ as the chain length $N$ grows. In fact,
$N=1$ already attains the displayed ratio exactly. After translation,
its point set is $\{0,\pm e_1,\ldots,\pm e_D\}$. Order the noncentral
vertices as
\[
e_1,e_2,\ldots,e_D,-e_1,-e_2,\ldots,-e_D.
\]
Consecutive vertices in this order are orthogonal, so their distances
are $\sqrt2$. Cut the order after its first and its second vertex,
creating three nonempty paths. Join the first vertex of each path to
the center. The resulting tree has maximum degree three, exactly three
unit edges, and $2D-3$ edges of length $\sqrt2$. Together with the
lower estimate, this proves
\[
\frac{\MSTthree(\{0,\pm e_1,\ldots,\pm e_D\})}
{\MST(\{0,\pm e_1,\ldots,\pm e_D\})}=L_D.
\]
Thus $L_D$ is an exact ratio for this particular point set, not a claim
that $\tau_3^D=L_D$. The chain supplies larger configurations with the
same guaranteed lower bound. In particular $L_D\to\sqrt2$ as $D\to\infty$.

\subsection{Orthogonal simplex clusters}
Choose integers $s\ge2$, $N\ge1$, and a real parameter $0\le\alpha<1$.
Here $s$ is the number of leaves in each cluster. Work in $sN$ dimensions
with coordinate unit vectors $f_{i,j}$, $1\le i\le N$, $1\le j\le s$.
Put
\[
u_i=\frac1{\sqrt s}\sum_{j=1}^s f_{i,j},\qquad
v_{i,j}=\sqrt s\,f_{i,j}-u_i.
\]
Then $\norm{u_i}=1$, $u_i\cdot v_{i,j}=0$, and
\[
\norm{v_{i,j}}^2=s-1,\qquad
v_{i,j}\cdot v_{i,k}=-1\quad(j\ne k),\qquad
\sum_{j=1}^s v_{i,j}=0.
\]
Thus the $v_{i,j}$ form a centered regular simplex in the
$(s-1)$-dimensional subspace perpendicular to $u_i$. Different cluster
subspaces are mutually orthogonal.

Take the root, centers, and leaves to be
\[
O=0,\qquad p_i=(1-\alpha)u_i,\qquad
q_{i,j}=u_i+v_{i,j}=\sqrt s\,f_{i,j}.
\]
There are $(s+1)N+1$ distinct points. For convenient distance notation set
\[
a=1-\alpha,\qquad A=\sqrt2\,a,\qquad
B=\sqrt{2s},\qquad H=\sqrt{s-1+\alpha^2}.
\]
All pair types have the following distances:
\begin{center}\renewcommand{\arraystretch}{1.1}
\begin{tabular}{lc}\toprule
Pair&Distance\\\midrule
$O,p_i$&$a$\\
$O,q_{i,j}$&$\sqrt s$\\
$p_i,p_h$ ($i\ne h$)&$A$\\
Two distinct leaves&$B$\\
$p_i,q_{i,j}$&$H$\\
$p_i,q_{h,j}$ ($i\ne h$)&$\sqrt{s+a^2}$\\\bottomrule
\end{tabular}\end{center}
For example, $u_i\cdot f_{i,j}=1/\sqrt s$, so
\[
\norm{p_i-q_{i,j}}^2
=a^2+s-2a=s-1+\alpha^2.
\]
All distinct leaf vectors are orthogonal and have squared norm $s$,
which gives their common distance $B$. For a center and another cluster's
leaf the dot product vanishes, giving squared distance $s+a^2$.
This exceeds $H^2$ by $2a\ge0$. The other entries follow directly from
orthogonality.

Joining $O$ to every $p_i$ and each $p_i$ to its $s$ leaves gives an
ordinary spanning tree of length
\begin{equation}\label{eq:clusterden}
N(a+sH).
\end{equation}
Only this upper bound on $\MST(X)$ is needed below.

\subsection{The degree-capacity charging argument}
We use Zbarsky's half-edge charging argument~\cite[Section~7]{Zbarsky14},
written here for a variable number $s$ of leaves.
Assume
\begin{equation}\label{eq:chargecondition}
A+B\le2H.
\end{equation}
Assign vertex charges
\[
c(O)=0,\qquad c(p_i)=A/2,\qquad c(q_{i,j})=B/2.
\]
The distance table shows that every distinct pair satisfies
\[
\norm{x-y}\ge c(x)+c(y).
\]
Center--center and leaf--leaf pairs give equality. Own center--leaf pairs
use \eqref{eq:chargecondition}; other center--leaf pairs are at least as
far apart. The root pairs use $a\ge A/2$ and $\sqrt s\ge B/2$.
Consequently every tree satisfies
\[
\len(T)\ge\sum_x\deg_T(x)c(x).
\]

A tree on these points has $(s+1)N$ edges, hence $2(s+1)N$ edge endpoints.
Let $t_0,t_p,t_q$ count the endpoints at the root, centers, and leaves.
For a degree-three tree,
\[
t_0\le3,\qquad t_p\le3N,\qquad
t_q=2(s+1)N-t_0-t_p.
\]
Since $B>A$, eliminating $t_q$ and then using these upper bounds gives
\begin{align*}
\len(T)&\ge\frac A2t_p+\frac B2t_q\\
&=\frac B2\bigl(2(s+1)N-t_0\bigr)+\frac{A-B}{2}t_p\\
&\ge\frac{3N}{2}A+\frac{(2s-1)N-3}{2}B.
\end{align*}
Combining this with \eqref{eq:clusterden} proves the finite bound
\begin{equation}\label{eq:Zfinite}
\tau_3^{sN}\ge Z_{s,N}(\alpha):=
\frac{\frac32 A+\frac12(2s-1-3/N)B}{a+sH}.
\end{equation}
The endpoint allocation is a relaxation of the degree constraints; it
need not be simultaneously attainable. Accordingly \eqref{eq:Zfinite}
is a lower bound, with no assertion that it is the exact ratio of the
cluster configuration.

\subsection{Zbarsky's choice and the limiting bound}
For integers $2\le s\le9$, choose
\begin{equation}\label{eq:alphas}
\alpha_s=-1-\sqrt s+\sqrt{4+4\sqrt s}.
\end{equation}
This lies in $[0,1)$. Indeed,
\[
4+4\sqrt s-(1+\sqrt s)^2=(3-\sqrt s)(1+\sqrt s)\ge0,
\]
whereas $4+4\sqrt s<(2+\sqrt s)^2$.
The equation $A+B=2H$, after division by $\sqrt2$ and squaring positive
sides, is equivalent to
\[
\alpha^2+2(1+\sqrt s)\alpha+s-3-2\sqrt s=0.
\]
Its root in $[0,1)$ is exactly \eqref{eq:alphas}. Thus the charge
condition holds with equality for every parameter used in the table below.

For $s=7$, these coordinates describe Zbarsky's seven-leaf clusters
\cite[Section~7]{Zbarsky14} in
an orthogonally rotated coordinate system. Each cluster consists of a
six-dimensional regular simplex and its perpendicular center direction,
so $N$ clusters require $7N$ dimensions. Its parameter is
\[
\alpha_7=-1-\sqrt7+\sqrt{4+4\sqrt7}
\approx0.173018810880.
\]
For fixed $s,\alpha$, define
\[
Z_s(\alpha)=\frac{\frac32\sqrt2(1-\alpha)
+\frac12(2s-1)\sqrt{2s}}
{1-\alpha+s\sqrt{s-1+\alpha^2}}.
\]
The finite formula has the exact correction
\begin{equation}\label{eq:Zcorrection}
Z_{s,N}(\alpha)=Z_s(\alpha)
-\frac{3\sqrt{2s}}{2N\bigl(1-\alpha+s\sqrt{s-1+\alpha^2}\bigr)}.
\end{equation}
Taking $N\to\infty$ therefore recovers Zbarsky's lower bound
\cite[Section~7]{Zbarsky14}
\[
\tau_3^\infty\ge Z_7(\alpha_7)
=1.447317759362\ldots.
\]
Precisely, for every $\varepsilon>0$, a sufficiently large finite $N$
makes the correction in \eqref{eq:Zcorrection} smaller than $\varepsilon$;
that finite configuration lives in dimension $7N$. This argument allows
the dimension to grow. In a prescribed dimension $D$, the constraint
$7N\le D$ requires the finite-$N$ formula instead.

\subsection{Comparison in a prescribed dimension}
Both families can be embedded in a larger ambient space by appending
zero coordinates. For $D\ge7$, the bound from Zbarsky's original choice~\cite[Section~7]{Zbarsky14} is
\[
\mathcal Z_7(D)=\max\{1,Z_{7,\lfloor D/7\rfloor}(\alpha_7)\}.
\]
For $D<7$ that particular family with $N\ge1$ does not fit, and its table
entry is a dash. To show the effect of changing the number of simplex
leaves, also define the explicit finite menu
\[
\mathcal Z(D)=\max\left(\{1\}\cup
\left\{Z_{s,\lfloor D/s\rfloor}(\alpha_s):
2\le s\le\min(9,D)\right\}\right).
\]
For fixed $s$ and $\alpha_s$, the finite bound increases with $N$, so
$N=\lfloor D/s\rfloor$ is the largest permitted, and best, choice of
$N$ for that fixed parameter. The maximum defining $\mathcal Z(D)$
checks only the displayed eight possible values of $s$; it does not
optimize over every construction or over the continuous parameter
$\alpha$. The automatic ratio bound one is included in both columns.

\begin{center}\small\renewcommand{\arraystretch}{1.08}
\begin{tabular}{rrrrrr}\toprule
Dimension $D$&Chain $L_D$&Zbarsky $\mathcal Z_7(D)$&Variable $s$: $\mathcal Z(D)$&$s$&$N$\\\midrule
2 & 1.1035533906 & --- & 1.0000000000 & --- & ---\\
3 & 1.2071067812 & --- & 1.0000000000 & --- & ---\\
4 & 1.2588834765 & --- & 1.0000000000 & --- & ---\\
7 & 1.3254535133 & 1.1357921403 & 1.1357921403 & 7 & 1\\
10 & 1.3520815280 & 1.1357921403 & 1.2233197492 & 5 & 2\\
14 & 1.3698335378 & 1.2915549498 & 1.2915549498 & 7 & 2\\
20 & 1.3831475452 & 1.2915549498 & 1.3331158444 & 5 & 4\\
28 & 1.3920235501 & 1.3694363546 & 1.3694363546 & 7 & 4\\
35 & 1.3964615526 & 1.3850126355 & 1.3850126355 & 7 & 5\\
40 & 1.3986805538 & 1.3850126355 & 1.3915114731 & 8 & 5\\
42 & 1.3994202209 & 1.3953968228 & 1.3953968228 & 7 & 6\\
47 & 1.4009939806 & 1.3953968228 & 1.3953968228 & 7 & 6\\
48 & 1.4012693885 & 1.3953968228 & 1.4013109576 & 6 & 8\\
49 & 1.4015335554 & 1.4028140995 & 1.4028140995 & 7 & 7\\
56 & 1.4031185562 & 1.4083770570 & 1.4083770570 & 7 & 8\\
70 & 1.4053375575 & 1.4161651975 & 1.4161651975 & 7 & 10\\
100 & 1.4080003589 & 1.4250659294 & 1.4250659294 & 7 & 14\\
200 & 1.4111069607 & 1.4361918444 & 1.4361918444 & 7 & 28\\
500 & 1.4129709217 & 1.4429300746 & 1.4429300746 & 7 & 71\\
1000 & 1.4135922420 & 1.4451239170 & 1.4451239170 & 7 & 142\\
\bottomrule\end{tabular}\end{center}
The entries are rounded for display; the exact lower bounds are the
formulas above. The final two columns give a maximizing parameter pair
within the stated finite menu. A dash there means that only the automatic
bound one is being displayed. The table compares these particular
families, not the best lower bounds known in each dimension.

The chain gives the stronger bound in the smaller listed dimensions.
At $D=48$, choosing $s=6,N=8$ gives $1.4013109576\ldots$, slightly above
$L_{48}=1.4012693885\ldots$. At $D=49$, Zbarsky's original choice
$s=7,N=7$ gives $1.4028140994\ldots$, above
$L_{49}=1.4015335553\ldots$. Larger dimensions permit more orthogonal
clusters. Ultimately the chain bound approaches $\sqrt2$, whereas the
seven-leaf cluster bound approaches $1.447317759362\ldots$.

Combining this known lower construction with the upper theorem gives
\[
Z_7(\alpha_7)\le\tau_3^\infty\le757/500.
\]
The separate rooted-path lower bound remains $3/2$. The interval for the
tree ratio and the interval for the rooted-path constant concern
different optimization problems.

\section*{Acknowledgement}
The authors acknowledge the use of ChatGPT (OpenAI) in exploring the
problem, developing and checking intermediate arguments, and preparing
the manuscript. The authors have reviewed the resulting arguments and
take full responsibility for the mathematical content and its presentation.

\begingroup\footnotesize\setlength{\parskip}{0pt}

\endgroup
\end{document}